\documentclass[12pt]{gtpart}
\usepackage{amsfonts, amsmath, amsthm}

\usepackage{graphicx}
\usepackage{epstopdf}
\usepackage{psfrag}

\title{Topological Index Theory for surfaces in 3-manifolds}

\author{David Bachman} 
\givenname{David}
\surname{Bachman}
\address{Pitzer College\\\newline
         Claremont, CA 91711\\USA}
\email{bachman@pitzer.edu}
\urladdr{http://pzacad.pitzer.edu/~dbachman}


\arxivreference{0901.0208}  
\arxivpassword{k5esy}   

\newtheorem{pro}{Proposition}[section]
\newtheorem{thm}[pro]{Theorem}
\newtheorem{lem}[pro]{Lemma}
\newtheorem{clm}[pro]{Claim}

\newtheorem{cnj}[pro]{Conjecture}
\newtheorem{cor}[pro]{Corollary}
\newtheorem{quest}[pro]{Question}

\theoremstyle{definition}
\newtheorem{dfn}[pro]{Definition}

\theoremstyle{remark}

\newcommand{\VV}{\mathcal V}
\newcommand{\WW}{\mathcal W}

\newcommand{\bdy}{\partial}

\begin{document}
\begin{abstract}
The disk complex of a surface in a 3-manifold is used to define its {\it topological index}. Surfaces with well-defined topological index are shown to generalize well known classes, such as incompressible, strongly irreducible, and critical surfaces. The main result is that one may always isotope a surface $H$ with topological index $n$ to meet an incompressible surface $F$ so that the sum of the indices of the components of $H \setminus N(F)$ is at most $n$. This theorem and its corollaries generalize many known results about surfaces in 3-manifolds, and often provides more efficient proofs. The paper concludes with a list of questions and conjectures, including a natural generalization of Hempel's {\it distance} to surfaces with topological index $\ge 2$. 
\end{abstract}
\maketitle

\section{Introduction.}

Let $H$ be a properly embedded, separating surface with no torus components in a compact, orientable 3-manifold $M$. Then the {\it disk complex},  $\Gamma(H)$, is defined as follows:
	\begin{enumerate}
		\item  Vertices of $\Gamma(H)$ are isotopy classes of compressions for $H$. 
		\item A set of $m+1$ vertices forms an $m$-simplex if there are representatives for each that are pairwise disjoint. 
	\end{enumerate}

Here we explore what information is contained in the topology of $\Gamma(H)$. To this end, we define:

\begin{dfn}
\label{d:Indexn}
The {\it homotopy index} of a complex $\Gamma$ is defined to be 0 if $\Gamma=\emptyset$, and the smallest $n$ such that $\pi_{n-1}(\Gamma)$ is non-trivial, otherwise. We say a surface $H$ is {\it topologically minimal} if its disk complex $\Gamma(H)$ is either empty or non-contractible. When $H$ is topologically minimal, we say its {\it topological index} is the homotopy index of $\Gamma(H)$. 
\end{dfn}

For example, a surface $H$ has topological index 1 if and only if $\pi_0(\Gamma(H))$ is non-trivial, i.e. its disk complex is disconnected. 

When $H$ is the boundary of a handlebody then the disk complex  was first defined by McCullough in \cite{mccullough:91}, who showed that in this case $\Gamma(H)$ is contractible. It follows that such surfaces are not topologically minimal. The goal of the present paper is to show that topologically minimal surfaces are a natural generalization of several well-known classes of surfaces in 3-manifolds, and that the results that hold for each of these classes also hold true for all topologically minimal surfaces. As an added benefit, proofs involving the set of all topologically minimal surfaces are often much shorter than existing proofs involving just, say, index 2 surfaces. This is largely owing to the inductive nature of the arguments. 

By definition, incompressible surfaces have topological index 0. In the next section we show that the {\it strongly irreducible} surfaces of Casson and Gordon \cite{cg:87} are precisely those that have topological index 1. We also show that {\it critical surfaces}, previously defined by the author in \cite{crit} and \cite{gordon}, have topological index 2. One important property shared by these types of surfaces is that they may always be isotoped to meet an incompressible surface in a collection of loops that are essential on both. We show here that this is in fact a corollary of a powerful result about all topologically minimal surfaces. This is given by Theorem \ref{c:MainCorollary}, which asserts that a topologically minimal surface $H$ and an incompressible surface $F$ can be isotoped so that $H \setminus N(F)$ is topologically minimal in $M \setminus N(F)$. 

Section \ref{s:Heegaard} contains corollaries to Theorem \ref{c:MainCorollary}. We show there that if $M$ contains a topologically minimal Heegaard surface then $\bdy M$ is incompressible. It then follows that if a closed 3-manifold $M$ contains any topologically minimal surface $H$ then either it is a Heegaard surface, $M$ is Haken, or $H$ is contained in a ball.  (In the final section we conjecture that this last possibility can not happen.) Finally, we show that if the disjoint union of surfaces is topologically minimal then so are its components, and its topological index is the sum of the indices of its components. Combining this with Theorem \ref{c:MainCorollary}, we find that a surface $H$ with topological index $n$ can be isotoped to meet an incompressible surface $F$ in such a way so that the sum of the indices of the components of $H \setminus N(F)$ is at most $n$. This is a generalization of known results about topological index 0 and 1 surfaces. 

In any new theory, the questions raised are as important as the new results. In the final section of this paper we list a few tantalizing questions and conjectures about topologically minimal surfaces. These include conjectures about the possible indices of topologically minimal surfaces in various kinds of 3-manifolds,  a natural generalization of Hempel's {\it distance} invariant \cite{hempel:01} to surfaces of arbitrary topological index, and conjectures which relate  geometric minimal surfaces to topologically minimal surfaces. 

Much of the motivation for this work comes from ideas of Hyam Rubinstein. In the late 1990's Rubinstein pioneered the viewpoint that strongly irreducible Heegaard splittings were the right class of surfaces within which to search for unstable (geometrically) minimal surfaces of index 1, as well as their PL analogues, the so-called ``almost normal" surfaces. One often finds such surfaces by minimax arguments involving 1-parameter sweepouts. Many of the topological arguments involving strongly irreducible surfaces also use 1-parameter sweepouts, so it became natural to think about such surfaces as being ``topologically minimal," in a very imprecise sense. In later work the author defined critical surfaces as an attempt to find some topological analogue to geometrically minimal surfaces that have index 2. As one would expect from such an analogue, arguments involving critical surfaces often involve 2-parameter sweepouts. In this paper we make precise the idea of topological index, demonstrate its usefulness, and conjecture its relation to  geometric minimal surfaces. 

The present work is the first in a sequence of papers on this topic. In \cite{barrier} we define a relative version of topological index for surfaces with non-empty boundary. The main result of that paper is that complicated amalgamating surfaces act as barriers to low index, low genus, topologically minimal surfaces. This is the key technical tool necessary for the author's  construction of a counter-example to the {\it Stabilization Conjecture} for Heegaard splittings \cite{StabilizationResults}. Further applications are given in \cite{AmalgamationResults}, where we prove several results about amalgamation and isotopy of Heegaard splittings. Finally, in joint work with Jesse Johnson, we produce examples of 3-manifolds containing surfaces with high topological index \cite{existence}. 

The author thanks several people for helpful comments during the preparation of this paper. Jesse Johnson had helpful suggestions regarding the construction of the family $H_x$ defined in the proof of Theorem \ref{t:MainTheorem}. Cameron Gordon, Daryl Cooper, and Andrew Casson provided advice necessary for the proof of Corollary \ref{c:IndexSum}. General  helpful comments were made by Martin Scharlemann and Yoav Moriah. Finally, it was Saul Schleimer and Eric Sedgiwck who first brought the index 1 case of Theorem \ref{c:ComponentsIntersection} to the attention of the author, which was the beginings of the paper \cite{glue}. In some sense this work is an extension of the main result of that paper.

\section{Low index surfaces}
\label{s:LowIndex}

In this section we show that the concept of topological index generalizes several well known classes of surfaces in 3-manifolds. 

\begin{dfn}
\label{d:H/D}
Let $H$ be a properly embedded surface in a 3-manifold $M$. A loop $\alpha$ on $H$ is {\it essential} if it does not bound a subdisk of $H$. A disk $D$ is a {\it compression} for $H$ if $D \cap H=\bdy D$ is an essential loop on $H$. The surface $H$ is {\it incompressible} if there are no compressions for it. If $D$ is a compression for $H$ then we construct the surface $H/D$ as follows. Let $M(H)$ denote the manifold obtained from $M$ by cutting open along $H$. Let $B$ denote a neighborhood of $D$ in $M(H)$. The surface $H/D$ is obtained from $H$ by removing $B \cap H$ and replacing it with the frontier of $B$ in $M(H)$. 
\end{dfn}

It follows immediately from the definitions that a surface has topological index 0 if and only if it is incompressible. We now show that surfaces with topological index 1 and 2 are also familiar.

Let $\VV$ and $\WW$ denote the sides of a Heegaard surface $H$, and $\Gamma_\VV(H)$ and $\Gamma_\WW(H)$ the subspaces of $\Gamma(H)$ spanned by compressions in $\VV$ and $\WW$. McCullough has called these complexes the {\it disk complexes} of $\VV$  and $\WW$. McCullough proved that such disk complexes are contractible  \cite{mccullough:91}. It follows that the topology of $\Gamma(H)$ is entirely determined by the simplices that connect $\Gamma_\VV(H)$ to $\Gamma_\WW(H)$. With this in mind, it is natural to introduce special terminology when there are no edges connecting $\Gamma_\VV(H)$ to $\Gamma_\WW(H)$. The following definition is due to Casson and Gordon \cite{cg:87}. 

\begin{dfn}
$H$ is {\it strongly irreducible} if there are compressions on opposite sides of $H$, but each compression on one side meets all compressions on the other. 
\end{dfn}

The main result of \cite{cg:87} is that if the minimal genus Heegaard splitting of a 3-manifold is not strongly irreducible, then the manifold contains an incompressible surface.

\begin{thm}
\label{l:TTStrongIrreducibility}
$H$ has topological index 1 if and only if it is strongly irreducible.
\end{thm}

\begin{proof}
By definition, a surface has topological index 1 when $\pi_0(\Gamma(H))$ is non-trivial. Hence, in this case $\Gamma(H)$ is disconnected. However, by McCullough's result $\Gamma_\VV(H)$ and $\Gamma_\WW(H)$ are contractible, so the only way for $\Gamma(H)$ to be disconnected is if both $\Gamma_\VV(H)$ and $\Gamma_\WW(H)$ are non-empty, and there are no edges connecting them. There are thus compressions on both sides, but any pair of such compressions intersect.
\end{proof}

In \cite{crit} the author introduced the idea of a {\it critical} surface. The main result of that paper is that if the minimal genus common stabilization of a pair of Heegaard splittings is not critical, then the manifold contains an incompressible surface. Critical surfaces were also instrumental in the author's proof of a conjecture of C. Gordon \cite{gordon}. 

\begin{dfn}
\label{d:critical}
$H$ is {\it critical} if the compressions for $H$ can be partitioned into sets $C_0$ and $C_1$ such that:
\begin{enumerate}
	\item For each $i=0,1$ there is at least one pair of disks $V_i, W_i \in C_i$ on opposite sides of $H$ such that $V_i \cap W_i =\emptyset$.
	\item If $V \in C_0$ and $W \in C_1$ are on opposite sides of $H$ then $V \cap W \ne \emptyset$.
\end{enumerate}
\end{dfn}

\begin{thm}
\label{l:critical}
$H$ has topological index 2 if and only if it is critical.
\end{thm}

\begin{proof}
We first establish that if $H$ has topological index 2 then it is critical.  Let $\Gamma_{\VV \WW}(H)$ be the subspace of $\Gamma(H)$ consisting of those cells spanned by vertices in both $\Gamma_\VV(H)$ and $\Gamma_\WW(H)$. 

\begin{clm}
\label{c:BetterPath}
Any path in $\Gamma_{\VV \WW}(H)$ which connects two vertices representing disks on the same side of $H$ is homotopic in $\Gamma_{\VV \WW}(H)$ to a path in either $\Gamma_\VV(H)$ or $\Gamma_\WW(H)$. 
\end{clm}

\begin{proof}
Let $\{D_i\}_{i=0}^n$ be such a path, where $D_0$ and $D_n$ are disks in $\VV$, and the path contains the fewest possible number of disks in $\WW$. Let $D_i$ be the first disk in this path in $\WW$. Then $D_{i-1}$ is in $\VV$. There are now two cases, depending on whether $D_{i+1}$ in $\VV$ or $\WW$. 

Consider first the case that $D_{i+1} \subset \VV$. We now produce a path $\{E^j\}_{j=0} ^m$ in $\Gamma_{\VV \WW}(H)$ from $D_{i-1}$ to $D_{i+1}$ that consists entirely of disks in $\VV$. By replacing $D_i$ with this path, we get a new path from $D_0$ to $D_n$ which contains fewer disks in $\WW$, a contradiction. 

First, note that by an innermost disk argument we may assume that each component of $D_{i-1}\cap D_{i+1}$ is an arc. Now we construct the new path from $D_{i-1}$ to $D_{i+1}$ by induction as follows:
	\begin{enumerate}
		\item Let $E^0=D_{i-1}$
		\item Assume we have constructed a sequence of disks $E^j$ in $\Gamma_{\VV \WW}(H)$ of disks that lie in $\VV$, such that for all $j$, $|E^{j} \cap D_{i+1}|<|E^{j-1} \cap D_{i+1}|$. If $E^{j} \cap D_{i+1}=\emptyset$ then we let $E^{j+1}=D_{i+1}$ and we have produced the desired path. Otherwise, let $\alpha$ be an arc of $E^{j} \cap D_{i+1}$ that is outermost on $D_{i+1}$. The arc $\alpha$ then cuts off a subdisk disk $D'_{i+1}$ of $D_{i+1}$ whose interior is disjoint from $E^{j}$. The arc $\alpha$ also cuts the disk $E^j$ into two subdisks. One of these, together with the disk $D'_{i+1}$, forms a compressing disk $E^{j+1}$ which meets $D_{i+1}$ fewer times than $E^j$ did. 
	\end{enumerate}
Note that every disk in the path $\{E^j\}$ is disjoint from $D_i$, a compressing disk in $\WW$. Hence, the entire path $\{E^j\}$ lies in $\Gamma_{\VV \WW}(H)$, as desired. 

We now move on to the case that $D_{i+1} \subset \WW$. As the edge $(D_i,D_{i+1})$ is in $\Gamma_{\VV \WW}(H)$, it must be part of a simplex $\Delta$ whose vertices represent disks on both sides of $H$. Therefore, there must be a vertex $E$ of $\Delta$ which represents a disk on the opposite side of $H$ as both $D_i$ and $D_{i+1}$, i.e. a disk in $\VV$. Now insert $E$ in the path between $D_i$ and $D_{i+1}$, and apply the above argument to again get rid of $D_i$. 
\end{proof}

Since $H$ has topological index 2 is follows immediately that $\pi_1(\Gamma(H))$ is non-trivial. Since both $\Gamma_\VV(H)$ and $\Gamma_\WW(H)$ are contractible we conclude there is a non-trivial loop $\gamma$ in $\Gamma(H)$ that passes from $\Gamma_\VV(H)$ to $\Gamma_\WW(H)$ and back, crossing through $\Gamma_{\VV \WW}(H)$ exactly twice. Let $(V_0,W_0)$ and $(V_1,W_1)$ be the two edges of $\Gamma_{\VV \WW}(H)$ traversed by this path, where $V_i \subset \VV$ and $W_i \subset \WW$. 

\begin{clm}
\label{c:GammaDisconnected}
The edges  $(V_0,W_0)$ and $(V_1,W_1)$ are in different components of $\Gamma_{\VV \WW}(H)$. 
\end{clm}

\begin{proof}
Suppose not. Then there is a path in $\Gamma_{\VV \WW}(H)$ connecting $V_0$ to $V_1$. By Claim \ref{c:BetterPath} there is such a path $\gamma'$ consisting entirely of disks in $\VV$. As $\Gamma_\VV(H)$ is simply connected, the path $\gamma \cap \Gamma _\VV (H)$ is homotopic to $\gamma'$. 

The edges $(V_0,W_0)$ and $(V_1,W_1)$, together with $\gamma'$, now form a path in $\Gamma_{\VV \WW}(H)$ from $W_0$ to $W_1$. By Claim \ref{c:BetterPath} this is homotpic to a path $\gamma''$ which lies entirely in $\Gamma _\WW(H)$. We have thus homotoped the original loop $\gamma$ entirely into $\Gamma_\WW(H)$, a simply connected space, contradicting its non-contractibility.
\end{proof}

Claim \ref{c:GammaDisconnected} immediately implies that $\Gamma_{\VV \WW}(H)$ is disconnected. We may therefore partition its components into two non-empty sets, $C_0$ and $C_1$, where $(V_i,W_i) \subset C_i$. Since $C_0$ and $C_1$ are a partition of the components of $\Gamma_{\VV \WW}(H)$, there are no edges $(V,W)$ that connect them, where $V \in C_0$ and $W \in C_1$. The sets $C_0$ and $C_1$ thus satisfy the conditions of Definition \ref{d:critical}. (Note that any vertex of $\Gamma(H)$ that is not in $\Gamma_{\VV \WW}(H)$ can be added to either $C_0$ or $C_1$, and the conditions of Definition \ref{d:critical} will still be satisfied.) 

\bigskip

To complete the proof of the theorem, we must now establish that if $H$ is critical then it has topological index 2. Let $C_i, V_i$, and $W_i$ be as in Definition \ref{d:critical}. We must produce a non-trivial loop in $\Gamma(H)$. Since $\Gamma_\VV(H)$ is contractible, there is a path of compressions in $\Gamma_\VV(H)$ from $V_0$ to $V_1$. Similarly, there is a path  from $W_0$ to $W_1$ in $\Gamma_\WW(H)$. These two paths, together with the edges $(V_i,W_i)$, form a loop $\alpha$ in $\Gamma(H)$. By way of contradiction, suppose $\alpha$ is trivial in $\pi_1(\Gamma(H))$. Then there is a map $f$ of a disk $D$ into $\Gamma(H)$ such that $f(\bdy D)=\alpha$. For some triangulation $T$ of $D$, we may assume $f$ is simplicial. We now assume that all choices have been made so that the number of 2-simplices in $T$ is minimal. 

Let $\Delta$ denote the triangle in $T$ that has $(V_0,W_0)$ as one of its edges. Without loss of generality we assume the third vertex of $\Delta$ represents a compression in $\VV$, and denote it as $V$. Since $(V,W_0)$ is an edge of $\Delta$, it follows that $V \cap W_0=\emptyset$. Hence, by criticality $V \in C_0$. If $V$ is in the interior of $D$ then remove $\Delta$ from $D$ and replace $V_0$ with $V$. This increases the combinatorial length of $\bdy D$, but reduces the number of 2-simplices in $T$, a contradiction. 

The remaining case is when $V$ is in $\bdy D$.  Then the edge $(V,W_0)$ cuts $D$ into two smaller disks. One of these, $D'$, contains the edge $(V_1,W_1)$. If we now replace $D$ with $D'$ and $V_0$ with $V$, we again contradict our minimality assumption.
\end{proof}

\section{Topological index in the complement of a surface}
\label{s:RelativeIndex}

In this section we show that a topologically minimal surface can always be isotoped so that it meets the complement of an incompressible surface in a topologically minimal surface. 

\begin{figure}
\psfrag{H}{$H$}
\psfrag{F}{$F$}
\psfrag{d}{$D'$}
\psfrag{D}{$D$}
\begin{center}
\includegraphics[width=3 in]{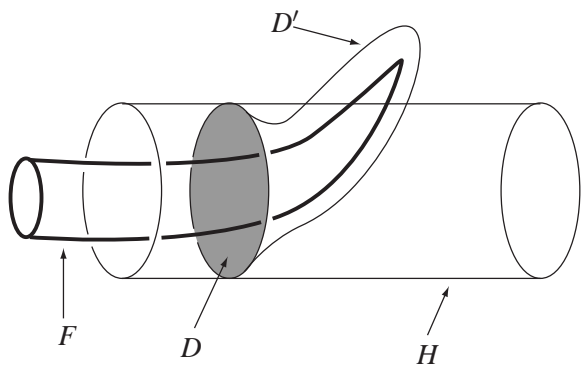}
\caption{A compression, $D$, for $H$ and its shadow, $D'$.}
\label{f:shadow}
\end{center}
\end{figure}

\begin{dfn}	
Let $H$ and $F$ be properly embedded surfaces in a 3-manifold $M$. Let $D$ be a compression for $H$. We say $D$ has a {\it shadow} (with respect to $F$) if there is a disk $D'$ where $\bdy D'=\bdy D$, $D' \cap F=\emptyset$, and the interior of $D'$ meets $H$ in loops  that are inessential on $H$.  The disk $D'$ is said to be a {\it shadow} of $D$. See Figure \ref{f:shadow}.
\end{dfn}

The main idea behind this paper is to exploit relationships between the homotopy indices of various complexes that depend on a specific position of $H$.  The first of these is the disk complex $\Gamma(H)$ of $H$. The complex $\Gamma_F(H)$ is the subset of $\Gamma(H)$ such that each vertex has a shadow. Later we will encounter a third complex, $\Gamma(H^F)$. 

The relationship between the homotopy indices of the complexes $\Gamma(H)$ and $\Gamma_F(H)$ is given presently in Theorem \ref{t:MainTheorem}. Later in this section we will use this theorem to prove that when $H$ is topologically minimal, then it can be isotoped so that it is topologically minimal in the complement of $F$. We then show that many of the standard results in 3-manifold topology, presently known for surfaces with low topological index, generalize to surfaces with arbitrary topological index. 



\begin{thm}
\label{t:MainTheorem}
Let $H$ and $F$ be properly embedded surfaces in $M$, where $H$ has topological index $n$. Then $H$ may be isotoped so that 
	\begin{enumerate}
		\item $H$ meets $F$ in $p$ points of tangency, for some $p \le n$. Away from these tangencies $H$ is transverse to $F$. 
		\item The complex $\Gamma_F(H)$ has homotopy index $i \le n-p$.
	\end{enumerate}
\end{thm}


\begin{proof}
When $H$ has topological index 0  the result is immediate, as $\Gamma_F(H) \subset \Gamma(H)=\emptyset$. We will assume, then, that $H$ has topological index $n \ge 1$. It follows that $\pi_{n-1}(\Gamma(H))$  is non-trivial, and thus there is a map $\iota \co S \to \Gamma(H)$ of an $(n-1)$-sphere $S$ into the $(n-1)$-skeleton of $\Gamma(H)$ which is not homotopic to a point. Let $B$ be the cone on $S$ to a point $z$. Hence, $B$ is an $n$-ball. 

Our first challenge is to define a continuous family of surfaces $H_x$ in $M$ isotopic to $H$, where $x \in B$.  Let $T$ be a triangulation of $S=\bdy B$ so that the map $\iota$ is simplicial. Let $\{v_i\}$ denote the set of vertices of $\Gamma(H)$ that are contained in $\iota(S)$. For each $i$ choose a representative $D_i$ from the isotopy class of disks represented by $v_i$ so that if $(v_i,v_j)$ is an edge of $\Gamma(H)$, then $D_i \cap D_j=\emptyset$. For each $i$, let $N_i$ be a small enough neighborhood of $D_i$ in $M$ so that $N_i \cap N_j=\emptyset$ whenever $(v_i,v_j)$ is an edge of $\Gamma(H)$, and let $f_i$ be a homeomorphism that takes $N_i$ to the standard unit ball in $\mathbb R^3=\{(x_1,x_2,x_3)\}$. Choose $f_i$ so that $f_i(H \cap N_i)$ is the graph of $r=1$ (in cylindrical coordinates), and $f_i(D_i)$ is a disk in the $x_1x_2$-plane. For each disk $D_i$ we now define a family of surfaces $H_i(t)$ in $N_i$, parameterized by a variable $t \in [0,1]$. These surfaces are given by the images of the graphs of $r=tx _3^2+1-t$, under the map $f_i^{-1}$ (see Figure \ref{f:H(t)}).

\begin{figure}
\begin{center}
\[\includegraphics[width=5 in]{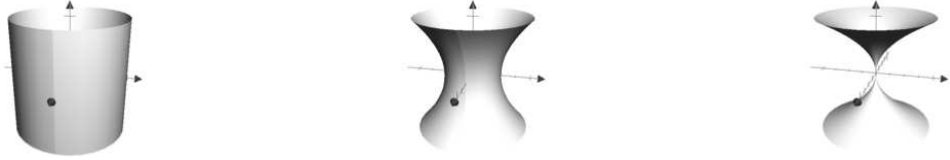}\]
\caption{The surfaces $H_i(t)$, for $t=0$, $\frac{1}{2}$, and 1.}
\label{f:H(t)}
\end{center}
\end{figure}

Extend $T$ to a triangulation $T'$ on $B$ by coning each simplex of $T$ to the point $z$. Suppose $\{D_0,...,D_{n-1}\}$ is the image of an $(n-1)$-simplex $\Delta$ of $T$ under the map $\iota$. We now identify the $n$-simplex of $T'$ which is the cone on $\Delta$ with the unit cube in $\mathbb R^n$. Label the axes of $\mathbb R^n$ with the variables $t_0, ..., t_{n-1}$. Place $z$ at the origin, and the vertex $v$ of $\Delta$ such that $\iota(v)=D_i$ at the point with $t_i=1$ and $t_j=0$ for all $j \ne i$. If $p$ is at the barycenter of a face $\sigma$ of $\Delta$ then place it at the vertex of the cube where the coordinates corresponding to the vertices of $\sigma$ are $1$ and the other coordinates are 0. We now linearly extend over the entire simplex to complete the identification with the cube. Now, if $x$ is in this $n$-simplex then $x$ has coordinates $(t_0(x),...,t_{n-1}(x))$. Let $H_x$ be the surface obtained from $H$ by replacing $H \cap N_i$ with the surface $H_i(t_i(x))$, for each $i$ between $0$ and $n-1$. See Figure \ref{f:HxDefn}. Repeating this for each $n$-simplex of $T'$ gives us the complete family of surfaces $H_x$. 

\begin{figure}
\psfrag{a}{$D_0$}
\psfrag{b}{$D_1$}
\psfrag{c}{$D_2$}
\psfrag{A}{$t_0$}
\psfrag{B}{$t_1$}
\psfrag{C}{$t_2$}
\psfrag{z}{$z$}
\psfrag{H}{$H$}
\psfrag{S}{$S$}
\begin{center}
\includegraphics[width=5 in]{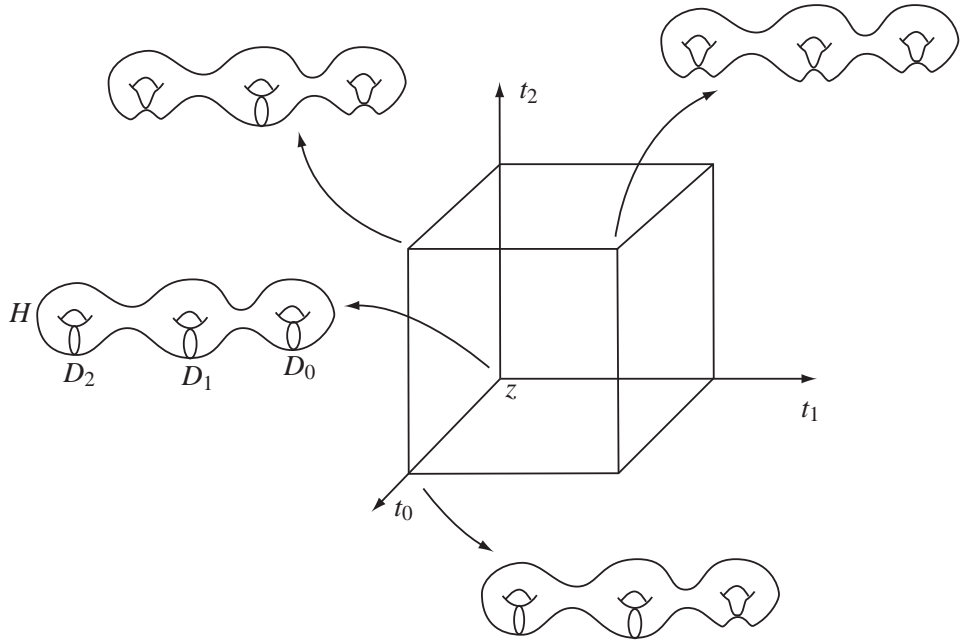}
\caption{A simplex $\Delta$ of $T'$, and a few of the surfaces $H_x$ for $x \in \Delta$. The union of the faces of the cube that do not meet $z$ is a simplex of $T$.}
\label{f:HxDefn}
\end{center}
\end{figure}

We assume $H$ is initially transverse to $F$. For each $i$, the surface $H_i(t) \subset N_i$ is tangent to $F$ for finitely many values $\{t_i^j\}$ of $t$. Hence, for each $x \in B$ the surface $H_x$ is tangent to $F$ at finitely many points, and each such point is in a distinct ball $N_i$. Note also that if $t_i(x)=t_i(y)$, then $H_x$ and $H_y$ agree inside of $N_i$. Hence, if $H_x$ is tangent to $F$ in $N_i$ then the surface $H_y$ will also be tangent to $F$, for all $y$ in the plane where $t_i(y)=t_i(x)$. It follows that each $n$-simplex of $T'$ is cubed by the points $x$ where $H_x$ is tangent to $F$. See Figure \ref{f:Cubing}. Hence, $B$ is cubed by the $n$-simplices of $T'$, together with this cubing of each such simplex. We denote this cubing of $B$ as $\Sigma$. It follows that if $x$ is in a codimension $p$ cell of $\Sigma$ then the surface $H_x$ is tangent to $F$ in at most $p$ points. 

\begin{figure}
\psfrag{A}{$t_0$}
\psfrag{B}{$t_1$}
\psfrag{C}{$t_2$}
\psfrag{t}{$\{t_0^j\}$}
\begin{center}
\includegraphics[width=3 in]{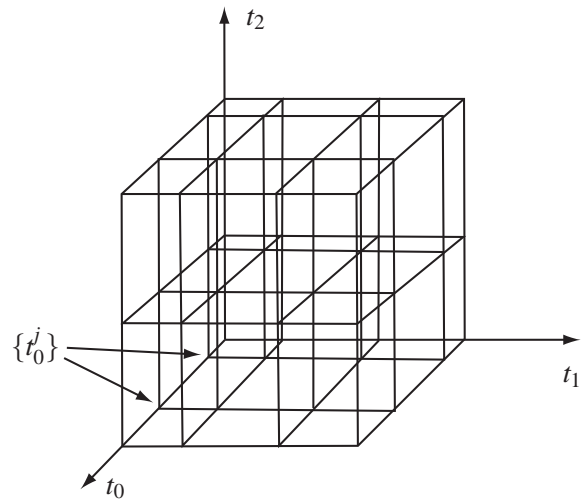}
\caption{A simplex $\Delta$ of $T'$ is cut up by planes into subcubes. Each such plane is determined by the points $x$ in which $H_x$ is tangent to $F$ in $N_i$, for some $i$.}
\label{f:Cubing}
\end{center}
\end{figure}

We now produce a contradiction by defining a continuous map $\Psi$ from $B$ into $\Gamma(H)$. The map $\Psi |\bdy B$ will be equal to $\iota$ on the barycenters of the $(n-1)$-cells of $T$, which will in turn imply that $\Psi$ maps $S$ onto $\iota(S)$ with the same degree as $\iota$. A contradiction follows as $\iota(S)$ is not homotopic to a point. 

For each $x \in B$ let $V_x=\Gamma_F(H_x)$. If $\tau$ is a cell of $\Sigma$, then we define $V_\tau$ to be the set $V_x$, for any choice of $x$ in the interior of $\tau$. Note that if $x$ and $y$ are in the interior of the same cell $\tau$ of $\Sigma$, then the pair $(H_x,F)$ is isotopic to $(H_y,F)$. Hence, $V_x=V_y$, and thus $V_\tau$ is well defined. The map $\Psi$ defined below will take each cell $\tau$ of $\Sigma$ into $V_\tau$. First, we establish a few properties of $V_\tau$. 

\begin{clm}
\label{c:subset}
Suppose $\sigma$ is a cell of $\Sigma$ which lies on the boundary of a cell $\tau$. Then $V_\sigma \subset V_\tau$. 
\end{clm}

\begin{proof}
Pick $x \in \sigma$ and $y \in \tau$. If $D \in V_x$ then $D$ is isotopic to a compression for $H_x$ that has a shadow $D'$.  To show $D \in V_y$ we must show that $D$ is isotopic to a compression for $H_y$ that has a shadow. Note that $H_y \cap F$ is obtained from $H_x \cap F$ by resolving some tangency. Hence, any loop of $H_x \setminus F$ is isotopic to a loop  of $H_y \setminus F$. It follows that since $\bdy D=\bdy D'$ was a loop  on $H_x$ disjoint from $F$, then $\bdy D=\bdy D'$ will be a loop  on $H_y$ that is disjoint from $F$. Furthermore, as the interior of $D'$ meets $H_x$ in a collection of loops  that are inessential on $H_x$, it follows that the interior of $D'$ meets $H_y$ in a collection of loops  that are inessential on $H_y$. We conclude $D'$ is a shadow for $D$, both as a compression for $H_x$ and as a compression for $H_y$. Hence, $D \in \Gamma_F(H_y)=V_y$. 
\end{proof}

\begin{clm}
\label{c:contradiction}
For each cell $\tau$ of $\Sigma$, 
\[\pi_{i}(V_\tau) = 1 \mbox{ for all }i \le {\rm dim}(\tau)-1.\]
\end{clm}

\begin{proof}
Let $x$ be in the interior of a codimension $p$ cell $\tau$ of $\Sigma$. Then the dimension ${\rm dim}(\tau)$  is $n-p$. The surface $H_x$ is tangent to $F$ in at most $p$ points, and is transverse to $H_x$ elsewhere. Recall $V_x=\Gamma_F(H_x)$. Thus, if the theorem is false then $V_x$  is non-empty, and $\pi_i(V_x)=1$ for all $i \le n-p-1={\rm dim}(\tau)-1$. 
\end{proof}

We now define $\Psi$ on the $0$-skeleton of $\Sigma$. For each $0$-cell $x \in \Sigma$, we will choose a point in $V_x$ to be $\Psi(x)$. If $x$ is in the interior of $B$ then $\Psi(x)$ may be chosen to be an arbitrary point of $V_x$. If $x$ is a point of $S=\bdy B$ then $x$ is contained in (perhaps more than one) $(n-1)$-simplex $\Delta_x$ of $T$. Let $\Delta'_x$ denote the face of $\Delta_x$ spanned by the vertices $v$ such that $t_i(v)=1$ if $t_i(x)=1$, and $t_i(v)=0$ otherwise. (Note that if $x$ was on the boundary of $\Delta_x$, so that it was also contained in some other $(n-1)$-simplex of $T$, then we still end up with the same simplex $\Delta_x'$ of $T$.) So, for example, if $x$ is at the barycenter of $\Delta_x$ then $\Delta '_x=\Delta _x$. By construction, for each vertex $v$ of $\Delta'_x$ the surface $H_x$ is pinched to a point along a disk $D$ in the isotopy class of $\iota(v)$. Hence, for all $y$ near $x$ the disk $D$ is a compression for $H_y$ that is disjoint from $F$. It follows that the entire simplex $\iota(\Delta'_x)$ is contained in $V_x$, and thus we may choose the barycenter of $\iota(\Delta'_x)$ to be the image of $\Psi(x)$. In particular, if $x$ is the barycenter of $\Delta_x$ then $\Psi(x)=\iota(x)$. 

We now proceed to define the rest of the map $\Psi$ by induction. Let $\tau$ be a $d$-dimensional cell of $\Sigma$. By induction,  assume $\Psi$ has been defined on the $(d-1)$-skeleton of $\Sigma$. In particular, $\Psi$ has been defined on $\bdy \tau$. Suppose $\sigma$ is a face of $\tau$.  By Claim \ref{c:subset} $V_\sigma \subset V_\tau$. By assumption $\Psi|\sigma$ is defined and $\Psi(\sigma) \subset V_\sigma$. We conclude $\Psi(\sigma) \subset V_\tau$ for all $\sigma \subset \bdy \tau$, and thus
\begin{equation}
\label{e:boundary}
\Psi(\bdy \tau) \subset V_\tau.
\end{equation}
Since $d={\rm dim}(\tau)$ it follows from Claim \ref{c:contradiction} that $\pi_{(d-1)}(V_\tau)=1$. Since $d-1$ is the dimension of $\bdy \tau$, we can thus extend $\Psi$ to a map from $\tau$ into $V_\tau$. 

What remains to be shown is that if $\tau$ is in $S=\bdy B$ then the extension of $\Psi$ from $\bdy \tau$ to $\tau$ may be done so that $\Psi(\tau) \subset \iota(S)$. Let $\Delta_\tau$ be the simplex of $T$ whose interior contains the interior of $\tau$. We need only show that $\Psi(\bdy \tau) \subset V_\tau \cap \iota(\Delta _\tau)$. Since $V_\tau \cap \iota(\Delta _\tau)$ will be a subsimplex of $\iota(\Delta_\tau)$, it follows that $\Psi$ can be extended over $\tau$ to this subsimplex. 

By Equation \ref{e:boundary}, $\Psi(\bdy \tau) \subset V_\tau$. So all we must do now is to show $\Psi(\bdy \tau) \subset \iota(\Delta _\tau)$. Let $\sigma$ denote a face of $\tau$, and $\Delta _\sigma$ the simplex of $T$ whose interior contains the interior of $\sigma$. Then $\Delta_\sigma$ is contained in $\Delta _\tau$. By induction we may assume $\Psi(\sigma) \subset \iota(\Delta _\sigma)$. Putting this together we conclude $\Psi(\sigma) \subset \iota(\Delta _\tau)$ for each $\sigma \subset \bdy \tau$, and thus $\Psi(\bdy \tau) \subset \iota(\Delta _\tau)$. 
\end{proof}

\begin{dfn}
Let $F$ be a properly embedded surface in a 3-manifold $M$.  Then we let $M^F$ denote the complement of a neighborhood of $F$ in $M$. For each subset $X$ of $M$, let $X^F=X \cap M^F$. 
\end{dfn}

\begin{figure}
\psfrag{A}{$\Gamma_F(H)$}
\psfrag{B}{$\Gamma(H^F)$}
\psfrag{C}{$\Gamma(H)$}
\psfrag{b}{Compressions for $H^F$ that are not compressions for $H$.}
\psfrag{c}{Compressions for $H$ contained in $M^F$.}
\begin{center}
\includegraphics[width=4 in]{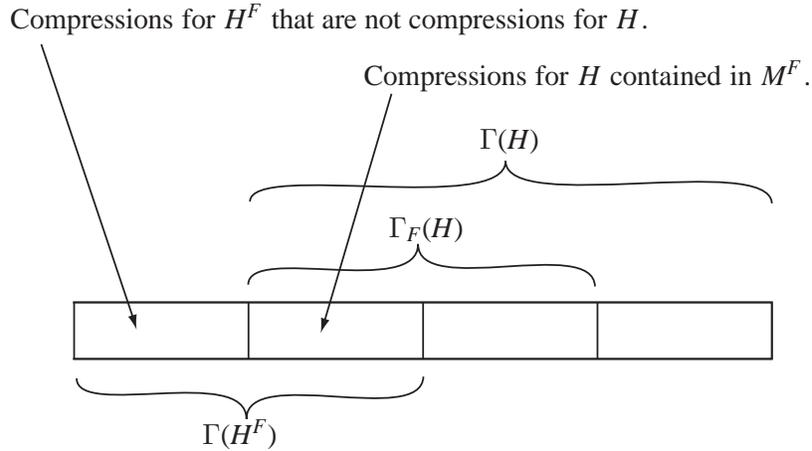}
\caption{Schematic showing how the complexes $\Gamma(H)$, $\Gamma_F(H)$, and $\Gamma(H^F)$ overlap.}
\label{f:complexes}
\end{center}
\end{figure}

We define the complex $\Gamma(H^F)$ precisely as above, where the vertices of $\Gamma(H^F)$ correspond to the compressions for $H^F$ in $M^F$. The relationship between the complexes $\Gamma(H)$, $\Gamma(H^F)$, and $\Gamma_F(H)$ is depicted in Figure \ref{f:complexes}.

We now use Theorem  \ref{t:MainTheorem} to show that when $H$ is topologically minimal and $F$ is incompressible, then $H$ may be isotoped so that $H^F$ is topologically minimal  in $M^F$. In Section \ref{s:Heegaard} we explore the implications of this when $H$ is a Heegaard surface.

\begin{lem}
\label{l:NonCompressionEffect}
Let $F$ be a properly embedded, incompressible surface in an irreducible 3-manifold $M$. Let $H$ be a properly embedded surface in $M$ which meets $F$ transversally, with the exception of a finite number of center and saddle tangencies, such that $\Gamma_F(H)$ has well defined homotopy index. Let $D$ be a compression for $H^F$ in $M^F$ that is not a compression for $H$. Then $\Gamma_F(H/D)=\Gamma_F(H)$.
\end{lem}

\begin{proof}
Let $M(H)$ and $B$ be as given in Definition \ref{d:H/D}. Then $H/D$ is obtained from $H$ by removing $B \cap H$ from $H$ and replacing it with the frontier $D_*$ of $B$ in $M(H)$. As $D$ is not a compression for $H$, $\bdy D$ bounds a subdisk $\mathcal D \subset H$. 

We first show  $\Gamma_F(H/D) \subset \Gamma_F(H)$. Suppose $E \in \Gamma_F(H/D)$. Then $\bdy E$ can be isotoped off of $D_*$. If $E$ now meets the ball $B$ then it can be further isotoped so that $E \cap B$ is a collection of disks parallel to $D$. But then each component of $E \cap B$ can be swapped with a disk parallel to $\mathcal D$. The resulting disk has the same boundary as $E$, but is disjoint from $H$. By the irreducibility of $M$ this disk must therefore be properly isotopic to $E$. See Figure \ref{f:HtoH/Dpart2}. We conclude that $E$ was a compression for $H$ that persisted as a compression for $H/D$. $E$ is therefore a compression for $H$ that is disjoint from $D$.

\begin{figure}
\psfrag{D}{$D$}
\psfrag{1}{$\mathcal D$}
\psfrag{d}{$D_*$}
\psfrag{H}{$H$}
\psfrag{h}{$H/D$}
\psfrag{F}{$F$}
\psfrag{E}{$E$}
\psfrag{e}{$E'$}
\begin{center}
\includegraphics[width=4.5 in]{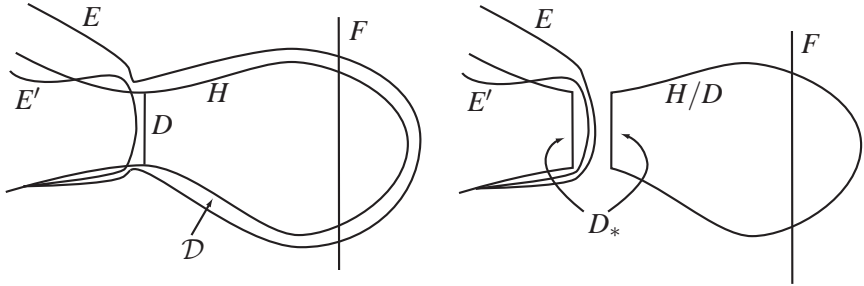}
\caption{Since $D$ is not a compression for $H$, any compression $E$ for $H/D$ (right figure) is always isotopic to a compression for $H$ (left figure). If $E'$ is a shadow for $E$ as a compression for $H/D$ (right figure), then $E'$ is a shadow for $E$ as a compression for $H$ (left figure).}
\label{f:HtoH/Dpart2}
\end{center}
\end{figure}

Now let $E'$ be a shadow for $E$ as a compression for $H/D$. As $\bdy E'=\bdy E$, it follows that $\bdy E' \cap D_*=\emptyset$. So, if $E'$ meets the ball $B$, then it meets it in disks parallel to $D$. The disk $E'$ thus meets $H$ in loops  isotopic to $E' \cap H/D$, together with loops  parallel to $D \cap H$. It follows that the interior of $E'$ meets $H$ in inessential loops, and thus $E'$ is a shadow for $E$ as a compression for $H$, i.e. $E \in \Gamma_F(H)$. See Figure \ref{f:HtoH/Dpart2}.

We now show $\Gamma_F(H) \subset \Gamma_F(H/D)$. Let $E$ now denote an element of $\Gamma_F(H)$. Thus, $\bdy E \cap F=\emptyset$. We assume $E$ has been chosen so that $|E \cap \mathcal D|$ is minimal. First we suppose $E \cap \mathcal D=\emptyset$. If the interior of $E$ meets $D$ then we may surger it off by a standard innermost disk argument. So in this case we may assume $E \cap D=\emptyset$. Since $E$ is a compression for $H$ but $D$ is not, it now follows that $E$ is a compression for $H/D$. Any shadow for $E$ as a compression for $H$ will be a shadow for $E$ as a compression for $H/D$, and thus $E \in \Gamma_F(H/D)$. 

Finally, we consider the case $E \cap \mathcal D \ne \emptyset$. Our goal is to isotope $E$ to a compression $E_0 \in \Gamma_F(H)$ such that $|E_0 \cap \mathcal D|<|E \cap \mathcal D|$, contradicting our minimality assumption. 

Let $\gamma$ denote an arc of $\bdy E \cap \mathcal D$ that is outermost on $\mathcal D$. Then $\gamma$ cuts a disk $\mathcal D'$ off of $\mathcal D$ whose interior does not meet $E$. We can use the disk $\mathcal D'$ to guide an isotopy of both $E$ and its shadow $E'$ to a compression $E_0$ for $H$ and a disk $E_0'$ with $\bdy E_0 =\bdy E_0'$. See Figure \ref{f:EtoE0}. Note that $|E_0 \cap \mathcal D|<|E \cap \mathcal D|$. If $\mathcal D' \cap F=\emptyset$, then it follows from the fact that $E'$ was a shadow of $E$ that $E_0'$ will be a shadow of $E_0$. Thus, $E_0 \in \Gamma_F(H)$ as desired.

\begin{figure}
\psfrag{D}{$\mathcal D$}
\psfrag{C}{$D$}
\psfrag{d}{$\mathcal D'$}
\psfrag{H}{$H$}
\psfrag{F}{$E_0$}
\psfrag{f}{$E'_0$}
\psfrag{E}{$E$}
\psfrag{e}{$E'$}
\begin{center}
\includegraphics[width=4.5 in]{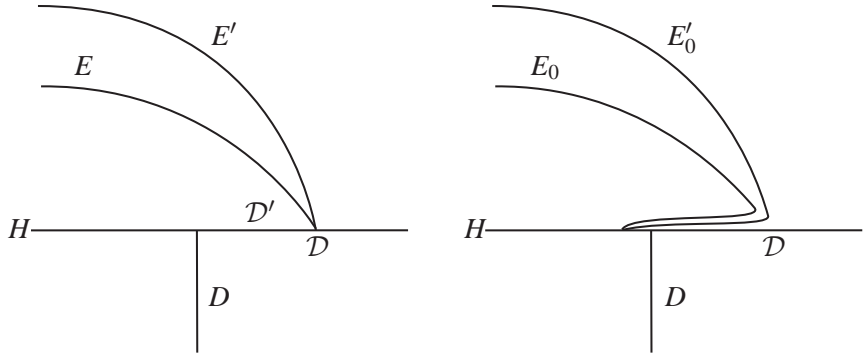}
\caption{Using the disk $\mathcal D'$ to obtain $E_0$ and $E_0'$ from $E$ and $E'$.}
\label{f:EtoE0}
\end{center}
\end{figure}

If $\mathcal D' \cap F \ne \emptyset$ then the disk $E_0'$ will not be a shadow for $E_0$, since $E_0' \cap F \ne \emptyset$. What remains then is to show that nonetheless, $E_0$ has a shadow. 

Let $N(F)$ denote a small product neighborhood of $F$. Since $E_0' \cap F \ne \emptyset$, it follows that $E_0' \cap \bdy N(F) \ne \emptyset$. Let $\delta$ denote a loop of $E_0' \cap \bdy N(F)$ that is outermost on $E_0'$. As $F$ is incompressible, $\delta$ bounds a subdisk $F_*$ of $\bdy N(F)$. See Figure \ref{f:E0toShadow}.

Although $F$ may not be transverse to $H$, the surface $\bdy N(F)$ will be. Thus, the disk $F_*$ meets $H$ in a collection of loops. We claim these loops are inessential on $H$, and thus $F_*$ can be used to surger $E_0'$ to a disk which meets $F$ fewer times. The new disk will meet $H$ more times, but each new intersection introduced will be inessential on $H$. Thus, by repeating this process we transform $E_0'$ to a shadow for $E_0$, as desired. 

To obtain a contradiction, suppose at least one loop of $F_* \cap H$ is essential on $H$. Let $\alpha$ be a such loop that is innermost on $F_*$. The loop $\alpha$ bounds a subdisk $A'$ of $F_*$ whose interior may meet $H$ in inessential loops. See Figure \ref{f:E0toShadow}. We claim $A'$ is the shadow of a compression $A$ for $H$, and thus $A \in \Gamma_F(H)$. 

\begin{figure}
\psfrag{f}{$E_0'$}
\psfrag{C}{$D$}
\psfrag{D}{$\mathcal D$}
\psfrag{H}{$H$}
\psfrag{G}{$F_*$}
\psfrag{B}{$B$}
\psfrag{A}{$A'$}
\psfrag{d}{$\delta$}
\psfrag{a}{$\alpha$}
\psfrag{b}{$\beta$}
\psfrag{F}{$N(F)$}
\begin{center}
\includegraphics[width=4.5 in]{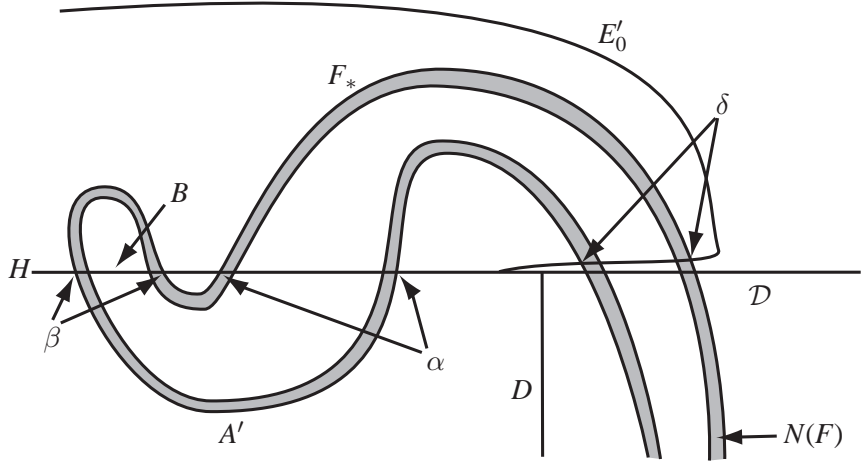}
\caption{The curves $\alpha$, $\beta$, and $\delta$, and the disks $A'$, $B$, and $F_*$.}
\label{f:E0toShadow}
\end{center}
\end{figure}

Let $\beta$ denote a loop of $A' \cap H$ that is innermost on $A'$. As $\beta$ is inessential on $H$ it bounds a subdisk $B$ of $H$. See Figure \ref{f:E0toShadow}. The disk $B$ can be used to surger $A'$, lowering $|A' \cap H|$. continuing in this way we arrive at a disk $A$ with the same boundary as $A'$ but whose interior is disjoint from $H$. As $\bdy A'=\bdy A$ is essential on $H$, we conclude $A$ is a compression for $H$. The disk $A'$ is then a shadow for $A$, and thus $A \in \Gamma_F(H)$.

Finally, suppose $X$ is any other element of $\Gamma_F(H)$. As $\bdy A =\bdy A' \subset F_* \subset \bdy N(F)$ and $\bdy X \cap F =\emptyset$, it follows that $\bdy X \cap \bdy A=\emptyset$. By a standard innermost disk argument (and the irreducibility of $M$) we may isotope $X$ to remove any intersections of its interior with the interior of $A$. Thus, we may assume $A \cap X=\emptyset$. The disk $X$ is therefore connected to the disk $A$ by an edge in $\Gamma_F(H)$. As this holds for all disks $X \in \Gamma_F(H)$, we conclude $\Gamma_F(H)$ is contractible to $A$. As $\Gamma_F(H)$ is not contractible, we have reached a contradiction. 
\end{proof}

\begin{thm}
\label{c:MainCorollary}
Let $F$ be a properly embedded, incompressible surface in an irreducible 3-manifold $M$. Let $H$ be a properly embedded surface in $M$ with topological index $n$. Then $H$ may be isotoped so that 
	\begin{enumerate}
		\item $H$ meets $F$ in $p$ saddle tangencies, for some $p \le n$. Away from these tangencies $H$ is transverse to $F$. 
		\item $H^F$ has topological index $i$, for some $i \le n-p$. 
	\end{enumerate}
\end{thm}

\begin{proof}
We begin by isotoping $H$ so as to satisfy the conclusion of Theorem \ref{t:MainTheorem}. Hence, we assume $H$ is tangent to $F$ in $p$ points, and the homotopy index of  $\Gamma_F(H)$ is at most $n-p$.

Let $D$ be a compression for $H^F$ that is not a compression for $H$. Then $\bdy D$ bounds a subdisk $\mathcal D$ of $H$. By Lemma \ref{l:NonCompressionEffect},  $\Gamma_F(H/D)=\Gamma_F(H)$. The surface $H/D$ contains a component $H'$ isotopic to $H$ (by the irreducibility of $M$), and a surface isotopic to $D \cup \mathcal D$. Note that as $\mathcal D \cap F \ne \emptyset$, $H'$ meets $F$ fewer times than $H$ did. Thus, we may repeat the above procedure only finitely many times. Note also that this procedure will remove all center tangencies of $H$ with $F$. We arrive at a surface $H_*$ isotopic to $H$ with $\Gamma_F(H_*)=\Gamma_F(H)$, such that every compression for $H_*^F$ is also a compression for $H_*$. As such compressions lie in the complement of $F$, they are their own shadows. Hence, such compressions are elements of $\Gamma_F(H_*)$. We conclude $\Gamma(H_*^F) \subset \Gamma_F(H_*)$. We claim the opposite inclusion is true as well, and thus $\Gamma(H_*^F) = \Gamma_F(H_*)$. 

Suppose now $E \in \Gamma_F(H_*)$. Let $E'$ be a shadow of $E$. Let $\beta$ be a loop of $E' \cap H_*$ that is innermost on $E'$. Then $\beta$ bounds subdisks $C \subset E'$ and  $C' \subset H_*$. If $C' \cap F \ne \emptyset$, then $C$ is a compression for $H_*^F$ that is not a compression for $H_*$, a contradiction. We conclude $C' \cap F=\emptyset$. Since $E' \cap F=\emptyset$ and $C \subset E'$, we conclude $C \cap F=\emptyset$. The sphere $C \cup C'$ thus bounds a ball in the complement of $F$ that we can use to guide an isotopy of $C$ to $C'$. (This may remove other components of $E' \cap C$ as well.) We thus transform the disk $E'$ to a disk $E''$ such that $\bdy E''=\bdy E$, $E'' \cap F=\emptyset$, and $|E'' \cap H_*|<|E' \cap H_*|$. Continuing in this way we arrive at a compression for $H_*$ with the same boundary as $E$, which is disjoint from $F$. Thus $E \in \Gamma(H_*^F)$.

We have now produced a surface $H_*$, isotopic to $H$, such that \[\Gamma(H_*^F)=\Gamma_F(H_*)=\Gamma_F(H).\] Thus, the homotopy index of $\Gamma(H_*^F)$ is equal to the homotopy index of $\Gamma_F(H)$. 
\end{proof}




\begin{cor}
\label{c:EssentialIntersection}
Let $F$ be a properly embedded, incompressible surface in an irreducible 3-manifold $M$. Let $H$ be a properly embedded surface in $M$ with topological index $n$. Then $H$ and $F$ may be isotoped so that any loop of $H \cap F$ is essential on both surfaces. 
\end{cor}

When $H$ is a Heegaard surface whose topological index is one this is a well-known result that has been used extensively in the literature. It first appears in the literature in the above form as Lemma 6 of \cite{schultens:00}, although it is implicit in Theorem 1 of \cite{kobayashi:88} and is directly implied by Theorem 1.1 of \cite{rs:97}. 

\begin{proof}
The first step is to use Theorem \ref{c:MainCorollary} to isotope $H$ so that $H^F$ is topologically minimal. The manifold $M^F$ is obtained from $M$ by removing a submanifold $N(F) \cong F \times I$. Let $F^1$ and $F^2$ denote the copies of $F$ on the boundary of $N(F)$. Each component  of $H \cap F^1$ is a loop or arc of $\bdy H^F$. Hence, we must show that every loop of $\bdy H^F$ that is inessential on $F^1$ is inessential on $H^F$. 

If there is a loop of $\bdy H^F$ that is inessential on $F^1$ then there is such a loop $\alpha$ that bounds a subdisk $C$ of $F^1$ whose interior is disjoint from $H^F$. If $\alpha$ is essential on $H^F$ then $C$ is a compression for $H^F$. Now suppose $D$ is some other element of $\Gamma(H^F)$. As $C \subset F^1$, the disks $C$ and $D$ can be made disjoint in $M^F$, and hence $(D,C)$ is an edge of $\Gamma(H^F)$. We conclude $C$ is connected by an edge to every other element of $\Gamma(H^F)$. It follows that $\Gamma(H^F)$ is contractible to $C$, a contradiction. 

We conclude that all loops of $H \cap F^1$ that are inessential on $F^1$ are also inessential on $H$. Any such loop thus bounds a disk component of $H^F$ that can be isotoped into $N(F)$, without affecting $\Gamma(H^F)$. By successively performing this operation we thus arrive at the desired position of $H$ with respect to $F^1$, a surface isotopic to $F$. 
\end{proof}

\section{Heegaard surfaces}
\label{s:Heegaard}

In this section we give some applications of topological index theory to  Heegaard splittings of 3-manifolds. We also show that the topological index of a surface is the sum of the topological indices of its components. 

\begin{lem}
\label{l:NegativeBdy}
Let $H$ be a properly embedded surface which separates $M$ into $\VV$ and $\WW$. Let $H_\VV$  be a surface obtained from $H$ by a sequence of compressions into $\VV$. Then $H_\VV$ is incompressible in the submanifold cobounded by $H$ and $H_\VV$. 
\end{lem}

\begin{proof}
Let $\{D_i\}$ denote the union of the compressions used to obtain $H_\VV$ from $H$. Let $E$ denote a compression for $H_\VV$ that lies between $H$ and $H_\VV$. By an innermost disk argument, we may surger $E$ off of each disk $D_i$. But the complement of a neighborhood of $\bigcup D_i$ in this submanifold is a product. As the boundary of a product does not admit compressions, we have thus reached a contradiction. 
\end{proof}

\begin{lem}
\label{l:TopMinSubmfld}
Let $H$ be a properly embedded surface which separates $M$ into $\VV$ and $\WW$. Let $H_\VV$ and $H_\WW$  be surfaces obtained from $H$ by maximal sequences of compressions into $\VV$ and $\WW$. Let $M_{\VV \WW}$ be the submanifold of $M$ cobounded by $H_\VV$ and $H_\WW$. If $H$ is topologically minimal  in $M$ then $H$ is topologically minimal  in $M_{\VV \WW}$.
\end{lem}

\begin{proof}
It suffices to show that every compression for $H$ in $M$ is isotopic to a compression in $M_{\VV \WW}$. Let $D$ be such a compression, and assume $D \subset \VV$. Isotope $D$ so that it meets $H_\VV$ minimally. If $D\cap H_\VV=\emptyset$, then the conclusion of the lemma follows. Hence, we assume there is a subdisk $D'$ of $D$, cut off by $H_\VV$, whose interior is disjoint from $H_\VV$. If $D' \cap H_\VV$ is not essential, then we contradict our assumption that $|D \cap H_\VV|$ is minimal. Hence, $D' \cap H_\VV$ is essential and we conclude  $D'$ is a compression for $H_\VV$. 

If $D'$ lies outside of $M_{\VV \WW}$ then we contradict the maximality of the sequence of compressions used to obtain $H_\VV$. But if $D'$ lies in $M_{\VV \WW}$ then it is in the submanifold cobounded by $H$ and $H_\VV$. This contradicts Lemma \ref{l:NegativeBdy}.
\end{proof}

\begin{thm}
\label{t:IncompressibleBoundary}
Let $H$ be a properly embedded surface which separates $M$ into $\VV$ and $\WW$. Let $H_\VV$  be a surface obtained from $H$ by a maximal sequence of compressions into $\VV$. If $H$ is topologically minimal  then $H_\VV$ is incompressible in $M$. 
\end{thm}

\begin{proof}
Let $H_\WW$ be the surface obtained from $H$ by a maximal sequence of compressions into $\WW$, and $M_{\VV \WW}$ the submanifold of $M$ cobounded by $H_\VV$ and $H_\WW$. By Lemma \ref{l:TopMinSubmfld} the surface $H$ is topologically minimal  in $M_{\VV \WW}$. 

We now claim that if either $H_\VV$ or $H_\WW$ is compressible, then there is a compression for one that misses the other. Assume there is no such compression for $H_\WW$. Let $D$ be a compression for $H_\VV$ in $M$. Isotope $D$ so that it meets $H_\WW$ minimally. If $D$ misses $H_\WW$ then we establish our claim. Assume then that $D$ meets $H_\WW$. Let $D'$ be a subdisk of $D$ cut off by $H_\WW$. If $\bdy D'$ is inessential on $H_\WW$, then we contradict our assumption that $|D \cap H_\WW|$ is minimal. But if $\bdy D'$ is essential on $H_\WW$ then $D'$ is a compression for $H_\WW$ that misses $H_\VV$, a contradiction. We conclude there is a compression $D$ for either $H_\VV$ or $H_\WW$ that misses the other. That is, $D$ is a compression for $H_\VV \cup H_\WW$. 

If $D$ lies outside of $M_{\VV \WW}$ then we contradict the minimality of the sequence of compressions used to obtain $H_\VV$ or $H_\WW$. Hence, $D \subset M_{\VV \WW}$. Note that $D$ is itself a properly embedded, incompressible surface in $M_{\VV \WW}$. We may thus apply Corollary \ref{c:EssentialIntersection} to isotope $H$ in $M_{\VV \WW}$ to meet $D$ in a collection of loops that are essential on both surfaces. Since $D$ does not contain any essential loops, we conclude $D \cap H=\emptyset$. 

The disk $D$ now lies either between $H$ and $H_\VV$, or between $H$ and $H_\WW$. In either case we contradict Lemma \ref{l:NegativeBdy}.
\end{proof}

\begin{cor}
Let $H$ be a topologically minimal Heegaard surface in a 3-manifold, $M$. Then $\bdy M$ is incompressible. 
\end{cor}

In the topological index one case this follows also from a celebrated Lemma of Haken \cite{haken:68}. In the topological index two case it was established by the author in \cite{gordon}. 

\begin{proof}
Let $\VV$, $\WW$, $H_\VV$, and $H_\WW$ be as in Theorem \ref{t:IncompressibleBoundary}. Since $H$ is a Heegaard surface, every component of $\bdy M$ is parallel to a component of either $H_\VV$ or $H_\WW$. The result is thus an immediate application of Theorem \ref{t:IncompressibleBoundary}. 
\end{proof}

\begin{cor}
\label{c:TopMinHeegaard}
Let $H$ be a closed topologically minimal surface in an irreducible  3-manifold, $M$. Then either 
	\begin{enumerate}
		\item $M$ contains a non-boundary parallel, incompressible surface, 
		\item $H$ is a Heegaard surface in $M$, 
		\item $H$ is contained in a ball, or 
		\item $H$ is isotopic into a neighborhood of $\bdy M$. 
	\end{enumerate}
\end{cor}

In the next section we conjecture that the third possibility does not happen. In particular, if $M$ is a non-Haken 3-manifold then it would follow that every topologically minimal surface in $M$ is a Heegaard surface. 

\begin{proof}
Let $\VV$, $\WW$, $H_\VV$, and $H_\WW$ be as in Theorem \ref{t:IncompressibleBoundary}. Suppose first some component of $H_\VV \cup H_\WW$ is a sphere. By the irreducibility of $M$, this sphere bounds a ball. If the ball contains $H$, then the result follows. Otherwise, we may remove each such sphere component from $H_\VV \cup H_\WW$. If the resulting surfaces are boundary parallel, then either $H$ is contained in a neighborhood of some boundary component of $M$, or $H$ is a Heegaard splitting of $M$. If some component of $H_\VV \cup H_\WW$ is not boundary parallel then by Theorem \ref{t:IncompressibleBoundary} it is incompressible, and the result follows. 
\end{proof}

\begin{lem}
\label{c:join}
Suppose $F$ and $G$ are disjoint surfaces in an irreducible 3-manifold $M$, and $F \cup G$ is topologically minimal. Then $\Gamma(F \cup G)$ is the join of $\Gamma(F)$ and $\Gamma(G)$. 
\end{lem}

\begin{proof}
Let $H=F \cup G$. Let $\VV$, $\WW$, $H_\VV$, and $H_\WW$ be as in Theorem \ref{t:IncompressibleBoundary}. By Theorem \ref{t:IncompressibleBoundary} the surfaces $H_\VV$ and $H_\WW$ are incompressible in $M$. 

If $E$ is a compression for $F$ then, as $H_\VV$ and $H_\WW$ are incompressible, we may isotope $E$ so that it is disjoint from both of these surfaces. It follows that $E$ is entirely contained in the component of $M_{\VV \WW}$ that contains $F$. But the surfaces $F$ and $G$ lie in different components of $M_{\VV \WW}$. Thus, $E$ must be disjoint from the surface $G$. Hence, any compression for $F$ is isotopic to a compression for $F \cup G$. We conclude there is a one-to-one correspondence between the vertices of $\Gamma(H)$ and the vertices of $\Gamma(F) \cup \Gamma(G)$. As every compression for $F$ will be disjoint from every compression for $G$, we conclude that $\Gamma(H)$ is the join of $\Gamma(F)$ and $\Gamma(G)$. 
\end{proof}

\begin{thm}
\label{c:IndexSum}
Suppose $F$ and $G$ are disjoint separating surfaces in an irreducible 3-manifold $M$, and $F \cup G$ is topologically minimal. Then $F$ and $G$ are topologically minimal  and
\[{\rm ind}(F)+{\rm ind}(G)={\rm ind}(F \cup G).\]
\end{thm}

Note that the hypothesis that $F \cup G$ is topologically minimal is extremely important. For example, let $F$ and $G$ be parallel surfaces in $M$ that each have topological index one. Then all of the compressing disks for $H=F \cup G$ are on the same ``side" of $H$. Hence, by McCullough's result \cite{mccullough:91}, $\Gamma(H \cup G)$ is contractible. Thus $H$ does not have topological index two, as one might expect.

\begin{proof}
We first show that $F$ and $G$ are topologically minimal. If not, then $\Gamma(F)$ (say) is non-empty and contractible. But the join of a contractible space with any other space is also contractible. It thus follows from Corollary \ref{c:join} that $F \cup G$ is not topologically minimal. 

If either $F$ or $G$ has topological index 0  then the result is immediate. We assume, then, that the topological index of $F$  is $n \ge 1$ and the topological index of $G$  is $m \ge 1$.

By definition, $(n-1)$ is the smallest $i$ such that $\pi_i(\Gamma(F)) \ne 1$, and $(m-1)$ is the smallest $j$ such that  $\pi_j(\Gamma(G)) \ne 1$. Our goal is to show that $(n+m-1)$ is the smallest $k$ such that $\pi_{k} (\Gamma(F \cup G)) \ne 1$. By Corollary \ref{c:join}, this is equivalent to showing that $(n+m-1)$ is the smallest $k$ such that $\pi_{k} (\Gamma(F) * \Gamma(G)) \ne 1$.

When $n=2$ then $\pi_1(\Gamma(F)) \ne 1$. Suppose $F$ seperates $M$ into $\VV$ and $\WW$. Let $\Gamma_\VV(F)$ and $\Gamma_\WW(F)$ denote the subsets of $\Gamma(F)$ spanned by the compressions that lie in $\VV$ and $\WW$, respectively. By an argument identical to the one given by McCullough in \cite{mccullough:91}, $\Gamma_\VV(F)$ and $\Gamma_\WW(F)$ are contractible. If we contract these to points $p_\VV$ and $p_\WW$, then the remaining 1-simplices of $\Gamma(F)$ join these two points. The fundamental group $\pi_1(\Gamma(F))$ is generated by these 1-simplices. The remaining 2-simplices have become bigons that run once over each of two 1-simplices. Hence, each such 2-simplex gives rise to a relation in $\pi_1(\Gamma(F))$ that kills one generator. It follows that $\pi_1(\Gamma(F))$ is free, and hence the non-triviality of $\pi_1(\Gamma(F))$ implies $H_1(\Gamma(F))$ is also non-trivial. Similarly, if $m=2$ we conclude $H_1(\Gamma(G))$ is non-trivial. For $n \ge 3$ the non-triviality of $H_{n-1}(\Gamma(F))$ follows from the Hurewicz Theorem. 

By Lemma 2.1 from \cite{milnor}:
\begin{eqnarray*}
\tilde H_{n+m-1} (\Gamma(F) &*& \Gamma(G))\\
 & \cong & \sum \limits _{i+j=n+m-2} \tilde H_i(\Gamma(F)) \otimes \tilde H_j(\Gamma(G))\\
&& + \sum \limits _{i+j=n+m-3} {\rm Tor}(\tilde H_i(\Gamma(F)), \tilde H_j(\Gamma(G))).
\end{eqnarray*}

In particular, it follows from the fact that $(n-1)$ is the smallest $i$ such that $H_i(\Gamma(F))$ is non-trivial, and $(m-1)$ is the smallest $j$ such that  $H_j(\Gamma(G))$ is non-trivial, that $(n+m-1)$ is the smallest $k$ such that $H_{k} (\Gamma(F) * \Gamma(G))$ is non-trivial. 
\end{proof}

As an immediate corollary we obtain:

\begin{cor}
\label{c:IndexNComponents}
If the topological index of $H$ is $n$, then the sum of the indices of the components of $H$ is exactly $n$. \qed
\end{cor}

Combining Theorem \ref{t:MainTheorem} with Corollary \ref{c:IndexNComponents} implies:

\begin{thm}
\label{c:ComponentsIntersection}
Let $F$ be a properly embedded, incompressible surface in an irreducible 3-manifold $M$. Let $H$ be a properly embedded surface in $M$ with topological index $n$. Then $H$ may be isotoped so that
	\begin{enumerate}
		\item $H$ meets $F$ in $p$ saddles, for some $p \le n$, and 
		\item the sum of the topological indices of the components of $H^F$, plus $p$, is at most $n$. 
	\end{enumerate}
\end{thm}

When $H$ is a Heegaard surface whose topological index is one, this result says that $F$ cuts $H$ up into incompressible pieces, along with at most one index one piece. Versions of this result were obtained by Schultens for graph manifolds \cite{schultens:04}, and the author, Sedgwick, and Schleimer for more general Haken manifolds \cite{glue}. 

Note also the similarity to the classification of {\it almost normal} surfaces given by Rubinstein. Such surfaces are cut up by the 2-skeleton of a triangulation into triangles and quadrilaterals, and exactly one ``special" piece. Rubinstein \cite{rubinstein:93} and Stocking \cite{stocking:96} proved that topological index 1 surfaces can always be isotoped to be almost normal. We believe the analogy is not a  coincidence; A relative (with respect to the 1-skeleton of a triangulation) version of Theorem \ref{c:ComponentsIntersection} should recover the Rubinstein-Stocking result, and generalize it to arbitrary topological index.

\section{Questions}
\label{s:Questions}

In any new theory, the questions raised are as important as the new results. Here we compile a list of questions and conjectures that we hope will stimulate further research on topologically minimal surfaces. 

\begin{quest}
How does topological index behave under finite covers? Are covers of topologically minimal surfaces also topologically minimal? 
\end{quest}

In \cite{existence} we produce a fairly generic family of 3-manifolds $\{M_n\}$ such that $M_n$ contains a Heegaard surface $H_n$ whose topological index is precisely $n$. The manifold $M_n$ is the $n$-fold cover of $M_1$, and the surface $H_n$ is the lift of $H_1$. 

\begin{quest}
Does every manifold have a topologically minimal Heegaard splitting?
\end{quest}

\begin{quest}
Are there non-Haken 3-manifolds with surfaces that have topological index $\ge 3$?
\end{quest}

\begin{cnj}
Suppose $M$ contains unstabilized Heegaard surfaces $F$ and $G$ that do not have topological index 1. Suppose further that the minimal genus common stabilization of $F$ and $G$ does not have topological index 2. Then $M$ contains a surface that has topological index $3$. 
\end{cnj}

By \cite{cg:87} such a manifold would be Haken, and so this conjecture compliments the question that precedes it. 

\begin{quest}
Is there a single 3-manifold that has surfaces of arbitrarily high topological index?
\end{quest}

\begin{cnj}
\label{cnj:S3}
$S^3$ and $B^3$ do not contain topologically minimal surfaces.
\end{cnj}

A corollary would be that handlebodies do not contain closed topologically minimal surfaces. Note also that this conjecture rules out the third conclusion given by Corollary \ref{c:TopMinHeegaard}.

\begin{cnj}
\label{cnj:product}
Let $F$ be a surface of positive genus. Then the only connected, topologically minimal surfaces in $F \times I$ are a single copy of $F$ and two copies of $F$ connected by an unknotted tube. 
\end{cnj}

By the argument given in the proof of Corollary \ref{c:TopMinHeegaard}, any topologically minimal surface in $F \times I$ would be a Heegaard surface or would be contained in a ball. The only incompressible (i.e. index 0) surface in $F \times I$ is a copy of $F$. By \cite{st:93}, the only strongly irreducible (i.e. index 1) Heegaard surface is two copies  of $F$ connected by an unknotted tube. So, if Conjecture \ref{cnj:S3} is true, then Conjecture \ref{cnj:product} is equivalent to the assertion that $F \times I$ contains no topologically minimal surfaces whose index is larger than one.

\begin{quest}
Does the conclusion of Corollary \ref{c:EssentialIntersection} hold if $F$ is topologically minimal, but not incompressible? 
\end{quest}

Rubinstein and Scharlemann have shown \cite{rs:96} that Corollary \ref{c:EssentialIntersection} holds when $H$ and $F$ both have topological index 1. This was instrumental in their proof that there is an upper bound on the smallest genus of a common stabilization of Heegaard surfaces $F$ and $G$, in terms of the genera of $F$ and $G$.

\begin{cnj}
If $H$ has topological index $n$ then it is isotopic to a geometrically minimal surface whose index is at most $n$. 
\end{cnj}

The index 0 case was proved by Freedman, Hass and Scott \cite{fhs:83}, and the index 1 case by Pitts and Rubinstein \cite{pr:87}. If true, it would indicate that topologically minimal surfaces are truly special. One would not expect, for example, a ``random" surface in a 3-manifold to be isotopic to a minimal surface. 






\begin{quest}
Suppose $H$ has topological index $n$. What information does ${\rm rank}(H_{n-1}(\Gamma(H))$ carry? What about other algebraic invariants of $\Gamma(H)$?
\end{quest}

\begin{quest}
[Generalized Hempel distance] For each surface $H$ there is a natural map of $\Gamma(H)$ into $\mathcal C(H)$, its curve complex, where the image of a compression $D$ is $\bdy D$. By \cite{harer:86}, $\mathcal C(H)$ has the homotopy type of a wedge of spheres. It follows that for low values of $n$ (in relation to the genus of $H$), each map $f \co S^{n-1} \to \Gamma(H)$ can be extended to a map $\hat f \co B^n \to \mathcal C(H)$. If we make all choices so that the number $d(n)$ of $n$-dimensional simplices in $\hat f(B^n)$ is minimal, then we get an interesting invariant when $f(S^{n-1})$ is not homotopic to a point in $\Gamma(H)$. When $H$ is a Heegaard surface that has topological index 1, Hempel called the invariant $d(1)$ the {\rm distance} of $H$ \cite{hempel:01}. Many interesting results have been obtained about Hempel's distance. What can be said about the invariant $d(n)$ for larger values of $n$?
\end{quest}

\bibliographystyle{alpha}

\begin{thebibliography}{McC91}

\bibitem[Baca]{barrier}
D.~Bachman.
\newblock Barriers to topologically minimal surfaces.
\newblock Available at {\tt http://arxiv.org/abs/0903.1692}.

\bibitem[Bacb]{StabilizationResults}
D.~Bachman.
\newblock Heegaard splittings of sufficiently complicated 3-manifolds {I}:
  {S}tabilization.
\newblock Available at {\tt http://arxiv.org/abs/0903.1695}.

\bibitem[Bacc]{AmalgamationResults}
D.~Bachman.
\newblock Heegaard splittings of sufficiently complicated 3-manifolds {II}:
  {A}malgamation.
\newblock Preprint.

\bibitem[Bac02]{crit}
D.~Bachman.
\newblock Critical {H}eegaard surfaces.
\newblock {\em Trans. Amer. Math. Soc.}, 354(10):4015--4042 (electronic), 2002.

\bibitem[Bac08]{gordon}
D.~Bachman.
\newblock Connected {S}ums of {U}nstabilized {H}eegaard {S}plittings are
  {U}nstabilized.
\newblock {\em Geometry \& Topology}, (12):2327--2378, 2008.

\bibitem[BJ]{existence}
D.~Bachman and J.~Johnson.
\newblock On the existence of high index topologically minimal surfaces.
\newblock In preparation.

\bibitem[BSS06]{glue}
David Bachman, Saul Schleimer, and Eric Sedgwick.
\newblock Sweepouts of amalgamated 3-manifolds.
\newblock {\em Algebr. Geom. Topol.}, 6:171--194 (electronic), 2006.

\bibitem[CG87]{cg:87}
A.~J. Casson and C.~McA. Gordon.
\newblock Reducing {H}eegaard splittings.
\newblock {\em Topology and its Applications}, 27:275--283, 1987.

\bibitem[FHS83]{fhs:83}
Michael Freedman, Joel Hass, and Peter Scott.
\newblock Least area incompressible surfaces in {$3$}-manifolds.
\newblock {\em Invent. Math.}, 71(3):609--642, 1983.

\bibitem[Hak68]{haken:68}
W.~Haken.
\newblock {\em Some Results on Surfaces in 3-Manifolds}.
\newblock M.A.A., Prentice Hall, 1968.

\bibitem[Har86]{harer:86}
John~L. Harer.
\newblock The virtual cohomological dimension of the mapping class group of an
  orientable surface.
\newblock {\em Invent. Math.}, 84(1):157--176, 1986.

\bibitem[Hem01]{hempel:01}
J.~Hempel.
\newblock 3-manifolds as viewed from the curve complex.
\newblock {\em Topology}, 40:631--657, 2001.

\bibitem[Kob88]{kobayashi:88}
Tsuyoshi Kobayashi.
\newblock Casson-{G}ordon's rectangle condition of {H}eegaard diagrams and
  incompressible tori in {$3$}-manifolds.
\newblock {\em Osaka J. Math.}, 25(3):553--573, 1988.

\bibitem[McC91]{mccullough:91}
Darryl McCullough.
\newblock Virtually geometrically finite mapping class groups of
  {$3$}-manifolds.
\newblock {\em J. Differential Geom.}, 33(1):1--65, 1991.

\bibitem[Mil68]{milnor}
J.~Milnor.
\newblock {\em Morse Theory}.
\newblock PUP, Princeton, New Jersey, 1968.

\bibitem[PR87]{pr:87}
J.~Pitts and J.~H. Rubinstein.
\newblock Applications of minimax to minimal surfaces and the topology of
  3-manifolds.
\newblock In {\em Miniconference on geometry and partial differential
  equations, 2 (Canberra 1986)}. Proc. ~Centre Math. ~Anal. ~Austral. ~Nat.
  ~Univ., {\bf 12}, Austral. ~Nat. ~Univ., Canberra, 1987.

\bibitem[RS96]{rs:96}
H.~Rubinstein and M.~Scharleman.
\newblock Comparing {H}eegaard {S}plittings of non-{H}aken 3-{M}anifolds.
\newblock {\em Topology}, 35:1005--1026, 1996.

\bibitem[RS97]{rs:97}
Hyam Rubinstein and Martin Scharlemann.
\newblock Transverse {H}eegaard splittings.
\newblock {\em Michigan Math. J.}, 44(1):69--83, 1997.

\bibitem[Rub95]{rubinstein:93}
J.~H. Rubinstein.
\newblock Polyhedral minimal surfaces, {H}eegaard splittings and decision
  problems for 3-dimensional manifolds.
\newblock In {\em Proceedings of the Georgia Topology Conference}, pages 1--20,
  1995.

\bibitem[Sch00]{schultens:00}
Jennifer Schultens.
\newblock Additivity of tunnel number for small knots.
\newblock {\em Comment. Math. Helv.}, 75(3):353--367, 2000.

\bibitem[Sch04]{schultens:04}
Jennifer Schultens.
\newblock Heegaard splittings of graph manifolds.
\newblock {\em Geom. Topol.}, 8:831--876 (electronic), 2004.

\bibitem[ST93]{st:93}
Martin Scharlemann and Abigail Thompson.
\newblock Heegaard splittings of {$({\rm surface})\times I$} are standard.
\newblock {\em Math. Ann.}, 295(3):549--564, 1993.

\bibitem[Sto00]{stocking:96}
M.~Stocking.
\newblock Almost normal surfaces in 3-manifolds.
\newblock {\em Trans. Amer. Math. Soc.}, 352:171--207, 2000.

\end{thebibliography}

\end{document}